\documentclass{article}
\usepackage{amsmath,amsfonts,amssymb,amsthm,pb-diagram,picinpar,graphicx,color}
\usepackage{algorithm}
\usepackage{amscd}
\usepackage{algorithmic}
\usepackage{colortbl}
\usepackage{enumerate}
\usepackage{hyperref}
\usepackage{bm}
\usepackage{tikz}
\usepackage{pgf}
\usepackage[T1]{fontenc}
\usepackage{enumitem}
\usepackage{comment}
\usepackage[hang,small,bf]{caption}
\usepackage[subrefformat=parens]{subcaption}
\usetikzlibrary{shapes,matrix,calc,arrows,decorations.markings,knots,patterns,intersections,cd}

\newcommand{\CC}{\mathbb{C}}%the field of complex numbers
\newcommand{\ZZ}{\mathbb{Z}}%the ring of integers
\newcommand{\PP}{\mathbb{P}}%projective space
\newcommand{\mcC}{\mathcal{C}}
\newcommand{\mcB}{\mathcal{B}}
\newcommand{\mcP}{\mathcal{P}}
\newcommand{\mcF}{\mathcal{F}}

\newcommand{\mcL}{\mathcal{L}}%tautological bundle

\newcommand{\mcO}{\mathcal{O}}%structure sheaf
\newcommand{\pic}{\mathop{\mathrm{Pic}}\nolimits}

\newcommand{\Supp}{\mathop{\mathrm{Supp}}\nolimits}

\newcommand{\sing}{\mathop{\mathrm{Sing}}\nolimits}

\newcommand{\comb}{\mathop{\mathrm{Comb}}\nolimits}

\newtheorem{thm}{Theorem}[section]
\newtheorem{lem}[thm]{Lemma}     % 2nd argument is what is printed

\newtheorem{prop}[thm]{Proposition}

\theoremstyle{definition}
\newtheorem{defin}[thm]{Definition}

\newtheorem{ex}[thm]{Example}

\theoremstyle{remark}
\newtheorem{rem}[thm]{Remark}

\title{Poncelet's closure theorem and the embedded topology of conic-line arrangements}
\author{Shinzo Bannai\footnote{Partially supported by JSPS KAKENHI Grant Numbers JP18K03263, JP23K03042} \and Ryosuke Masuya \and Taketo Shirane\footnote{Partially supported by JSPS KAKENHI Grant Number JP21K03182} \and Hiro-o Tokunaga\footnote{Partially supported by JSPS KAKENHI Grant Number JP20K03561} \and Emiko Yorisaki}
\date{\today}

\begin{document}

\maketitle

\abstract{In this paper, we consider conic-line arrangements that arise from Poncelet's closure theorem. We study unramified double covers of the union of two conics, that are induced by a $2m$-sided Poncelet transverse. As an application, we show the existence of  families of  Zariski pairs of degree $2m+6$ for $m\geq 2$ that consist of reducible curves having two conics and $2m+2$ lines as irreducible components. }

\section{Introduction}
There is a very famous theorem in projective geometry known as Poncelet's closure theorem or Poncelet's porism,  first introduced in \cite{poncelet}.
Let  $C_1, C_2\subset \PP^2$ be smooth conics that intersect transversely. Given a general point $P_1\in C_1$, let $L_1$ be a line passing through $P_1$ and tangent to $C_2$. Then $L_1$ intersects $C_1$ at another point $P_2$ and we can choose  a line $L_2$ passing through $P_2$ and tangent to $C_2$. By repeating this process, we obtain a sequence $(P_1, L_1), (P_2, L_2), \ldots$ of pairs of points $P_i\in C_1$ and tangent lines $L_i \in C_2^\ast$ such that $P_i\in L_i$, where $C_2^\ast$ is the dual curve of $C_2$. Such a sequence is called a \emph{Poncelet transverse originating at $P_1$}. Poncelet's closure theorem can be  stated as follows:
\begin{thm}[{\cite{poncelet}}]\label{thm:poncelet}
Let $C_1$, $C_2$ be as above. If there exists a point $P\in C_1$ such that the Poncelet transverse originating at $P$ is periodic with period $n$, then for any $Q\in C_1$, the Poncelet transverse originating at $Q\in C_1$ is also periodic with period $n$.  
\end{thm}
A modern proof in terms of torsion points of elliptic curves was given by P.~Griffiths and J.~Harris in \cite{griffiths-harris78}. We call a Poncelet transverse $\sqcap_n$ with minimum period $n$ an \emph{$n$-sided Poncelet transverse} and denote it by a sequence $\sqcap_n:(P_1, L_1), \ldots, (P_n, L_n)$ of $n$ pairs. Sometimes, we omit \lq$n$-sided Poncelet' and simply use \lq transverse' when it is evident from the context that we are talking about an $n$-sided Poncelet transverse. For a general choice of origin $P_1\in C_1$, an $n$-sided Poncelet transverse will consist of  $n$ distinct points $P_1, \ldots, P_n$  and $n$ distinct lines $L_1, \ldots, L_n$ and can be viewed as a $n$-gon with these points and lines as vertices and edges that is inscribed in $C_1$ and circumscribed about $C_2$. Poncelet's theorem also holds for special choices of origins, where  $\sqcap_n$ %the $n$-sided Poncelet transverse 
will involve  intersection points of $C_1$ and $C_2$ and/or bitangent lines of $C_1$ and $C_2$ and can be viewed as a degenerated $n$-gon having double edges.
In this paper, we consider  curves arising from $2m$-sided Poncelet transverses ($m\geq 2$) and construct a new infinite sequence of curves that are interesting with regard to their \emph{embedded topology}.
%We hope this  topological viewpoint  is  new  and will add to the already rich literature of curves realeted to Poncelet's closure theorem.    

The embedded topology of an algebraic plane curve $\mcC\subset \PP^2$ is the homeomorphism class of the pair $(\PP^2, \mcC)$ of topological spaces.  One of the main objectives is to give a complete classification of the embedded topology of plane curves. It is known that for  two plane curves $\mcC_1, \mcC_2\subset \PP^2$, if $(\PP^2, \mcC_1)$ and $(\PP^2, \mcC_2)$ are homeomorphic as pairs, then $\mcC_1$ and $\mcC_2$ have the same combinatorial type.
However the converse is not true, and there exist pairs of plane curves $\mcC_1, \mcC_2$  that have the same combinatorial type but the homeomorphism classes of  $(\PP^2, \mcC_1)$ and $(\PP^2, \mcC_2)$ are distinct. Such pairs  $(\mcC_1, \mcC_2)$ of curves are called \emph{Zariski pairs} and have been the interest of many mathematicians (see \cite{zariski29}, \cite{survey} for precise definitions and details). Understanding what causes the differences in the embedded topology of Zariski pairs should lead to a better understanding of the embedded topology of plane curves in general and hence is important. 
Concerning Zariski pairs of reducible plane curves with irreducible components of small degree, the following results are known. In the case where $\mcC_i$ are line arrangements, it is known that Zariski pairs do not exist for $\deg \mcC_i\leq 9$ (see \cite{Nazir-Yoshinaga2012}). Also, it is known that there exists a Zariski pair of line arrangements of $\deg\mcC_i=11$ (see \cite{Artal-Carmona-Cogolludo-Marco-05}). However the case of $\deg\mcC_i=10$ is open. In the case of conic-line arrangements, i.e. reducible curves whose irreducible components are lines and smooth conics, the existence of a Zariski pair of degree $7$ consisting of curves with two conics and three lines as irreducible components is known (see \cite{tokunaga14}). Also, a Zariski pair of degree $6$ consisting of curves with three conics as irreducible components is known (see \cite{Oka2007}). It would be interesting to find a Zariski pair of  conic-line arrangements of low degree consisting of curves with a single conic and additional lines.

Now, we explain our main result. Let $C_1$, $C_2$ be smooth conics intersecting transversely that admit a $2m$-sided Poncelet transverse. Let $L_1, \ldots, L_{2m}$ be distinct  lines forming a non-degenerated $2m$-sided Poncelet transverse $\sqcap_{2m}: (P_1, L_1), \ldots, (P_{2m}, L_{2m})$ where $L_1, \ldots, L_{2m}$ are tangent lines of $C_2$ and  the points $P_{1}=L_{2n}\cap L_1$ and $P_i=L_{i-1}\cap L_{i}$ ($2\leqq i \leqq 2n$)  lie on $C_1$. Let $\mcP_{2m}:=\sum_{i=1}^{2m} L_i$ be the union of the lines and let $T_1, T_2, T_3, T_4$ be the four bitangent lines of $C_1$ and $C_2$. The curves  that we are interested in are reducible curves of the form 
\[
\mcC_{ij}:=C_1+C_2+\mcP_{2m}+T_i+T_j \quad (\{i, j\}\subset \{1, 2, 3, 4\})
\]
having two conics and $2m+2$ lines as irreducible components. 
We consider the double covers of $\PP^2$ branched along the $2m+2$ lines $\mcP_{2m}+T_i+T_j$ and see how $C_1+C_2$ behaves under these double covers. Namely, we calculate the \emph{splitting types} of $C_1+C_2$ with respect to these double covers. (See Section~\ref{sec:splitting-type} and  \cite{bannai16} for details on splitting types.) This will be done by studying the invertible sheaves $\mcF$ of order 2 on $C_1+C_2$, or equivalently, torsion points of order 2 of the Jacobian $J(C_1+C_2)$ of the singular curve $C_1+C_2$.
As an application of these calculations,  we obtain the following theorem.
\begin{thm}\label{thm:main}
Under the above notation, for a suitable choice of labels of $T_1, \ldots, T_4$, the pair $(\mcC_{ij}, \mcC_{kl})$ is a Zariski pair if $\{i, j\}=\{1,2\}$ or $\{3, 4\}$ and $\{k, l\}\not=\{1,2\}, \{3, 4\}$.
\end{thm}
We note that in the above setting,  $C_1, C_2$ admits two degenerated $2m$-sided Pocelet transverses each having  two bitangent lines as edges. The differences of the curves $\mcC_{ij}$ and $\mcC_{kl}$ in the theorem are  whether the two bitangent lines lie in the same degenerated Poncelet transverse or not. 
Although the curves that are proved to be Zariski pairs in Theorem~\ref{thm:main} are conic-line arrangements of $\deg \mcC_{ij}\geq 10$, we believe that our method of systematically constructing Zariski pairs from Poncelet transverses is in itself interesting and worth sharing. We hope this  topological viewpoint  is  new  and will add to the already rich literature of curves realeted to Poncelet's closure theorem.

Similar studies relating torsion elements of the Jacobian $J(C)$ and the embedded topology of reducible curves having $C$ as an irreducible component have been done in \cite{bannai-tokunaga20-1},  \cite{BBST23-1}, \cite{BBST23-2} when $C$ is a smooth curve. This paper can be considered as a variation of these works in the case where $C$ is reducible and singular.

This paper is organized as follows: In Section~\ref{sec:splitting-type}, we review the definition of splitting types and state the proposition that is used in distinguishing the embedded topology. In Section~\ref{sec:linebdl2coverings}, we give a discussion on unramified double covers of conic-line arrangements, especially in the case of two transversal conics. In Section~\ref{sec:poncelet-cover}, we study $n$-sided Poncelet transverses and double covers related to them. Finally in Section~\ref{sec:proof-main}, we give the proof of our main result, Theorem~\ref{thm:main}.

%{\bf Acknowledgements:}

\section{Splitting types 
}\label{sec:splitting-type}
In this section, we review the notion of \emph{splitting types} of plane algebraic curves with respect to a double cover, which will be used to distinguish the embedded topology of the curves that we are interested in. We refer the reader to \cite{bannai16} for details.  Let $\pi_\mcB:S'\rightarrow \PP^2$ be a double cover branched along a curve $\mcB\subset \PP^2$ of even degree and let $C\subset\PP^2$ be an irreducible plane curve. The preimage $\pi_{\mcB}^{-1}(C)$ can be either reducible or irreducible, depending on the relation between $C$ and the branch locus $\mcB$. In the former case where $\pi_\mcB^{-1}(C)$ is reducible, $\pi_{\mcB}^{-1}(C)$ will have two irreducible components since $\pi_{\mcB}^{-1}$ is a double cover. In this case we say that  \emph{$C$ is a splitting curve with respect to $\pi_{\mcB}$ or $\mcB$}. Let $C_1, C_2$ be splitting curves with respect to $\mcB$ and let $\pi_{\mcB}^{-1}(C_i)=C_i^++C_i^-$, $(i=1, 2)$. The relation between the components $C_1^\pm$, $C_2^\pm$ reflect how the curves $\mcB, C_1, C_2$ are \lq intertwined' in $\PP^2$, hence gives information about the embedded topology of the reducible curve $\mcC=\mcB+C_1+C_2$. The information can be formulated as follows:
\begin{defin}\label{def:splitting-type}
Let $\mcB, C_1, C_2$ be as above. 
For integers $m_1\leq m_2$, we say that  the  triple $(C_1, C_2; \mcB)$ has splitting type $(m_1, m_2)$, if $C_1^+\cdot C_2^+=m_1$ and $C_1^+\cdot C_2^-=m_2$ for a suitable choice of labels.    
\end{defin}
The splitting types can be used to distinguish the embedded topology of reducible plane curves by the following proposition.
\begin{prop}[{\cite[Proposition~2.5]{bannai16}}]\label{prop:splitting-type}

Let $\mcB_1$, $\mcB_2$ be plane curves of degree $2d$ and let $C_{i1}, C_{i2}$ be splitting curves with respect to $\mcB_i$, $(i=1, 2)$. Suppose that $C_{i1}\cap C_{i2}\cap \mcB_i=\emptyset$, $C_{i1}$ and $C_{i2}$ intersect transversely and that $(C_{11}, C_{12};\mcB_1)$ and $(C_{21}, C_{22}; \mcB_2)$ have distinct spitting types. Then a homeomorphism $h:\PP^2\rightarrow \PP^2$ such that $h(\mcB_1)=\mcB_2$ and $\{h(C_{11}), h(C_{12})\}=\{C_
{21}, C_{22}\}$ does not exist.
\end{prop}

\begin{rem}
It is known that Definition~\ref{def:splitting-type} and Proposition~\ref{prop:splitting-type} can be modified to a more general version (see \cite{bannai-shirane-tokunaga-jaca}), but the above version is enough for our purposes. 
\end{rem}

Later, we will calculate the splitting types of pairs of conics $C_1, C_2$ with respect to various double covers in order to prove our main theorem.

\section{Line bundles of order two  and unramified double covers of conic-line arrangements with simple nodes}\label{sec:linebdl2coverings}

In this section we briefly recall the theory of double covers. We will especially consider the case of conic-line arrangements with simple nodes for later use. We refer the reader to \cite[Section~2, Section~3]{harris82} for details and arguments in a more general setting. 
% 記号の案
% \mcF : two torsion invertible sheaf on \mcC
% \mcL : invertible sheaf on \PP^2 s.t. \mcL^2=\mcO(\mcB)
% p_\mcL: \bm{L}_\mcL \to \PP^2 : the line bundle corresponding to \mcL

First, we consider topological (unramified) double covers. Let $\mcC=C_1+\dots+C_k$ be a conic-line arrangement with simple nodes, i.e. each irreducible component $C_i$ of $\mcC$ is either a line or a smooth conic, and all intersection points are ordinary double points.  Let $\varphi: \mcC'\rightarrow \mcC$ be a topological double cover of $\mcC$. Then, since $C_i\cong \PP^1$ and is simply connected, $\varphi^{-1}(C_i)$ splits into two disjoint sets  $C_i^\pm\subset \mcC'$. We fix a labeling $C_i^\pm$ for the meantime.  Let $Q\in \sing(\mcC)$,  and let $C_i, C_j$ be the irreducible components intersecting at $Q$. Then $C_i^+$ will intersect with either $C_j^+$ or $C_j^-$ over $Q$. We say that $\mcC'$ is glued by $+$ over $Q$ if $C_i^+$ intersects $C_j^+$ and is glued by $-$ over $Q$ if $C_i^+$ intersects $C_j^-$. We summarize this data in the form of a  map defined as below.
\begin{defin}
A \emph{gluing data of order two on $\mcC$} is  a map $\kappa:\sing(\mcC)\to\{+,-\}$. If there is no confusion, we simply call it a \textit{gluing data}. The \textit{gluing data $\kappa_{\varphi}$ of a topological double cover $\varphi:\mcC'\rightarrow \mcC$} is a gluing data on $\mcC$ defined by  $\kappa_{\varphi}(Q)=+$ if  $\mcC'$ is glued by $+$ over $Q$ and $\kappa_{\varphi}(Q)=-$ if  $\mcC'$ is glued by $-$ over $Q$.
\end{defin}
If we reverse the labeling of $C_i^\pm$ of  $\pi^{-1}(C_i)$, then all of the signs for $Q\in \sing(\mcC)\cap C_i$ will be reversed. Namely, for a gluing data $\kappa:\sing(\mcC)\to\{+,-\}$ and each $i=1,\dots,k$,  a new gluing data  $\kappa_{i}$  is obtained by 
\begin{align}\label{eq:gluing_data_operation} \kappa_{i}(Q):=\left\{ \begin{array}{ll}
    -\kappa(Q) & \mbox{if $Q\in C_i$,} \\
    \kappa(Q) & \mbox{otherwise}
\end{array} \right. \end{align}
for each $Q\in\sing(\mcC)$. We say that two gluing data $\kappa$ and $\kappa'$ are \textit{equivalent}, and write $\kappa\sim\kappa'$, if $\kappa'$ can be constructed from $\kappa$ by a finite number of the above operations. In this way, we have a map from the set of  topological double covers  $\varphi: \mcC'\rightarrow \mcC$ to the set of equivalence classes of gluing data $\kappa_{\varphi}$ on $\mcC$.

\begin{lem}\label{lem:rigidity-of-top-2cover}
Let $\mcC=C_1+\dots+C_k$ be a conic-line arrangement with simple nodes. The following map $\Psi$ from the set of homeomorphism class of topological double covers of $\mcC$ to the set of equivalence class of gluing data on $\mcC$ is well-defined and one-to-one: 
\[ \begin{array}{rccc}\Psi:& \{ \varphi : \mbox{a topological double cover of $\mcC$}\}/\cong & \to & \{\kappa: \mbox{a gluing data on $\mcC$}\}/\sim \\ & \rotatebox{90}{$\in$} & & \rotatebox{90}{$\in$} \\ & \lbrack \varphi\rbrack & \mapsto & \lbrack\kappa_\varphi\rbrack \end{array},  \]
where $\kappa_\varphi$ is the gluing data of the topological double cover $\varphi$ of $\mcC$. 
%
% $\varphi\mapsto\kappa_\varphi$ for each topological double cover $\varphi:\mcC'\to\mcC$ gives a one-to-one correspondence between the set of homeomorphism classes of topological double covers of $\mcC$ and the set of equivalence class of gluing data on $\mcC$. 
In particular, any continuous deformation $\varphi_t:\mcC'_t\to\mcC$ ($t\in\Delta$) of topological double covers of $\mcC$ is constant, where $\Delta\subset\CC$ is a small neighborhood of the origin. 
\end{lem}

\begin{proof}
% Let $\varphi:\tilde\mcC\to\mcC$ be a topological double cover. 
% Since $C_i\cong\PP^1$ is simply connected, $\varphi^{-1}(C_i)$ consists of two disjoint components $C_i^\pm$ for each $i=1,\dots,k$. 
% We define a gluing data $\kappa_\varphi$ on $\mcC$ by 
% \[ \kappa_\varphi(Q_{ij}^l):=\left\{ \begin{array}{ll}
%     + & \mbox{if $C_i^+$ and $C_j^+$ intersect over $Q_{ij}^l$,} \\
%     - & \mbox{otherwise}
% \end{array} \right. \]
% for each $i,j=1,\dots,k$ and $Q_{ij}^l\in C_i\cap C_j$. Note that $C_i\cap C_j$ consists of multiple points if $C_i$ or $C_j$ os a conic.
% For each $i=2,\dots,k$, fix a point $Q_{1i}^0\in C_1\cap C_i$. 

Suppose that $h:\mcC'\to\mcC''$ is a homeomorphism over $\mcC$ of topological double covers $\varphi':\mcC'\to\mcC$ and $\varphi'':\mcC''\to\mcC$. 
Put ${\varphi'}^{-1}(C_i)={C'}_i^++{C'}_i^-$ and ${\varphi''}^{-1}(C_i)={C_i''}^++{C_i''}^-$. 
Note that $h$ satisfies $h({\varphi'}^{-1}(C_i))={\varphi''}^{-1}(C_i)$. 
Since $h$ is a homeomorphism, the gluing data $\kappa_{\varphi''}$ is equivalent to $\kappa_{\varphi'}$ as $\kappa_{\varphi''}$ can be obtained from $\kappa_{\varphi'}$ by applying the operations (\ref{eq:gluing_data_operation}) to $\kappa_{\varphi'}$ for $1\leq i\leq k$ with $h({C'}_i^+)\ne {C_i''}^+$.
Hence the map $\Psi$ is well-defined. 

Let $\kappa:\sing(\mcC)\to\{+,-\}$ be a gluing data on $\mcC$. 
For each $i=1,\dots,k$, let $C_i^+\sqcup C_i^-$ be the disjoint union of two copies $C_i^\pm$ of $C_i$, and let $\varphi_i:C_i^+\sqcup C_i^-\to C_i$ be the projection, which is the topological double cover of $C_i$. 
We construct a topological space $\mcC'_\kappa$ by gluing $C_i^\pm\cap \varphi_i^{-1}(Q)$ to $C_j^\pm\cap\varphi_j^{-1}(Q)$ if $\kappa(Q)=+$, and to $C_j^\mp\cap\varphi_j^{-1}(Q)$ if $\kappa(Q)=-$, for each $1\leq i<j\leq k$ and $Q\in C_i\cap C_j$. 
Then the topological double covers $\varphi_i$ induce a topological double cover $\varphi_\kappa:{\mcC'}_\kappa\to\mcC$. 
Since the operation (\ref{eq:gluing_data_operation}) corresponds to replacing $C_i^+$ and $C_i^-$, the map $\lbrack\kappa\rbrack\mapsto \lbrack\varphi_\kappa\rbrack$ is well-defined, and is the inverse map of $\Psi$. 
Hence $\Psi$ is one-to-one. 
%Moreover, the maps $\varphi\mapsto\kappa_\varphi$ and $\kappa\mapsto\varphi_\kappa$ gives an expected correspondence since $\varphi_{\kappa_\varphi}=\varphi$ and $\kappa_{\varphi_\kappa}=\kappa$. 

Let $\varphi_t:\mcC'_t\to\mcC$ ($t\in\Delta$) be a continuous deformation of topological double covers, and let $\Phi:\overline{\mcC}'\to\Delta\times\mcC$ be the continuous family of the topological double covers $\varphi_t$, where $\overline{\mcC}':=\{(t, P')\mid t\in\Delta,\ {P'}\in{\mcC}'_t \}$ and $\Phi(t,P'):=(t,\varphi_t({P'}))$. 
Then $\Phi$ is a topological double cover of $\Delta\times \mcC$. 
Since $\Delta\times C_i$ is simply connected for each irreducible component $C_i\subset\mcC$, the preimage $\Phi^{-1}(\Delta\times C_i)$ consists of two connected components $\overline{C}_i^\pm$. 
For each $Q\in\sing(\mcC)$, the preimage $\Phi^{-1}(\Delta\times\{Q\})$ also consists of two components $\Delta_Q^\pm$. 
For each $t\in\Delta$, we define a gluing data $\kappa_t:\sing(\mcC)\to\{+,-\}$ by, for each $Q\in C_i\cap C_j$ ($i\ne j$), $\kappa_t(Q)=+$ if $\Delta_Q^+\subset \overline{C}_i^+\cap \overline{C}_j^+$, and $\kappa_t(Q)=-$ otherwise. 
This $\kappa_t$ coincides with the gluing data $\kappa_{\varphi_t}$ of $\varphi_t$ for any $t\in\Delta$. Since $\kappa_t$ is constant under $t$, all topological double covers $\varphi_t$ are homeomorphic. 
\end{proof}

Next, we consider the relation between unramified double covers of $\mcC$ and  invertible sheaves of order $2$ on $\mcC$ following \cite{harris82}. Let $\mcC$ be a reduced curve and let  $\mcF$ be an invertible sheaf of order $2$ on $\mcC$, i.e. $\mcF\otimes\mcF \cong \mcO_{\mcC}$ where $\mcO_\mcC$ is the structure sheaf of $\mcC$. Let $p_\mcF:\bm{L}_{\mcF}\rightarrow \mcC$ be the line bundle corresponding to  $\mcF$ and let $t\in \Gamma(\bm{L}_\mcF, p_\mcF^\ast \mcF)$ be the tautological section. Then, since we have assumed that  
$\mcF\otimes\mcF \cong \mcO_\mcC$, the zero divisor of $t^2-1$ in $\bm{L}_\mcF$ gives an  unramified double cover $\varpi_\mcF:\mcC'_\mcF\rightarrow \mcC$ of $\mcC$. Note that the construction is algebraic, but since it is unramified, $\varpi_\mcF$ is also a topological double cover. It is known that this relation induces a one-to-one correspondence between isomorphism classes of invertible sheaves $\mcF$ of  order 2, and isomorphism classes of unramified double covers $\varpi_\mcF:\mcC'_\mcF\rightarrow \mcC$ of $\mcC$. In the case where $\mcC$ is a conic-line arrangement with simple nodes, the  relation between invertible sheaves of order 2 and unramified double covers can be described using the gluing data as follows:  Let $\mcF$ be an invertible sheaf of order 2 on $\mcC$. Then, since the components of $\mcC$ are isomorphic to $\PP^1$ and $\pic(\PP^1)$ is torsion free, the restrictions $\mcF|_{C_i}$ are isomorphic to the trivial sheaf $\mcO_{C_i}$. Given, a node $Q\in C_i\cap C_j$, the line bundle $\mcF$ and its transition functions give the data of gluing of the  trivial sheaves $\mcF|_{C_i}$ and $\mcF|_{C_j}$ over $Q$ which is given by multiplication with $1$ or $-1$, since $\mcF$ is of order 2. Conversely, an invertible sheaf of order 2 can be constructed by assigning this gluing data $\pm1$ of the trivial sheaves at each node $Q$. In this way, we can associate a gluing data $\kappa_\mcF$ to $\mcF$. Again,  changing the signs of the gluing data at every node in an irreducible component $C_i$ will result in an isomorphic sheaf. Hence, we have the following Lemma:
\begin{lem}\label{lem:rigidity_of_sheaves} 
Let $\mcC=C_1+\cdots+C_k$ be a conic-line arrangement with simple nodes. Then, there is a one-to-one correspondence between the set of  invertible sheaves $\mcF$ of order $2$ on $\mcC$ and the set of equivalence classes of gluing data on $\mcC$. Furthermore, this correspondence is compatible with the correspondence between the unramified  double cover $\varpi_\mcF: \mcC'_\mcF\rightarrow \mcC$ of $\mcC$ associated to $\mcF$ and its gluing data. Namely if $\kappa_\mcF$ is the gluing data of $\mcF$ and $\kappa_{\varpi_\mcF}$ is the gluing data of $\varpi_\mcF: \mcC'_\mcF\rightarrow \mcC$, then $\kappa_\mcF \sim \kappa_{\varpi_\mcF}$. 
\end{lem}
\begin{proof}
The first part follows from the discussion stated before the lemma. See also \cite[Section~2b]{harris82}. For the second part, let $\mcF$ be an invertible sheaf of order 2 and let $\mcC'_\mcF\rightarrow \mcC$ be the associated unramified double cover. The preimage of $C_i$ in $\mcC'_\mcF$ consists of disjoint copies $C_i^+$ and $C_i^-$ of $C_i$ corresponding to the decomposition $t_i^2-1=(t_i-1)(t_i+1)$, where $t_i$ is the tautological section of $p^\ast(\mcF|_{C_i})$.
The  gluing data of $\mcF$ at a node $Q\in C_i\cap C_j$ tells us how the tautological sections $t_i, t_j$ are related and in turn how  $C_i^\pm$ and $C_j^\pm$ intersect as curves in $\bm{L}_\mcF$. If the sheaves $\mcF|_{C_i}$ and $\mcF|_{C_j}$  are glued by multiplication by $1$ over $Q$ then $t_i=t_j$ over $Q$  and $C_i^+$ intersects $C_j^+$ over $Q$. If they are glued by $-1$ then $t_i=-t_j$ over $Q$ and $C_i^+$ intersects $C_j^-$ over $Q$. Hence, the gluing data $\kappa_\mcF$ of  $\mcF$ as a sheaf and the gluing data  $\kappa_{\varpi_\mcF}$ of the topological unramified double cover $\varpi_\mcF: \mcC'_\mcF \to\mcC$ associated to $\mcF$ coincide. 
\end{proof}
Furthermore, given two invertible sheaves $\mcF_1, \mcF_2$ that are each of order $2$, the product $\mcF_1\otimes\mcF_2$ is again of order $2$. The gluing data of   $\mcF_1\otimes\mcF_2$  is given by simply taking the products of the gluing data of $\mcF_1$ and $\mcF_2$ at each $Q\in\sing(\mcC)$, as the transition functions of $\mcF_1\otimes\mcF_2$ are given by products of the transition function of $\mcF_1$ and $\mcF_2$. The gluing data of the  unramified double cover associated to $\mcF_1\otimes\mcF_2$ can also be obtained likewise.

Understanding and calculating the structure of the unramified double covers through this gluing data is useful and will be used in the proof of the main theorem.  Also, in some  cases where $\mcC\subset \PP^2$, and the unramified double cover of $\mcC$ is induced by a (possibly ramified) double cover of $\PP^2$, the structure of the former can be deduced from the latter as follows:
%この辺りcoveringとcoverが混在していたので統一
Let $\pi:S'\to\PP^2$ be a double cover branched along a plane curve $\mcB\subset\PP^2$ of degree $2d$, and let $F\in\Gamma(\PP^2,\mcO_{\PP^2}(2d))$ be a defining polynomial of $\mcB$. 
Assume that there is an effective divisor $D$ on the curve $\mcC$ such that $\mcB|_{\mcC}=2D$ and $\Supp D\cap\sing \mcC=\emptyset$, i.e., $\mcC$ intersects with $\mcB$ at smooth points of $\mcC$ with even multiplicities. 
Put $\mcL:=\mcO_{\PP^2}(d)$, and let $p_\mcL:\bm{L}_\mcL\to\PP^2$ be the line bundle corresponding to $\mcL$, where $\bm{L}_\mcL:=\mathrm{Spec}\,S(\mcL^{-1})$ is the spectrum of the symmetric algebra $S(\mcL^{-1})$ of $\mcL^{-1}$. 
Let $t\in \Gamma(\bm{L}_\mcL,p_\mcL^\ast\mcL)$ be the tautological section. 
Then $S'$ can be regarded as the subvariety of $\bm{L}_\mcL$ defined by $t^2-F=0$, and $\pi=p_\mcL|_{S'}$. 
Since $\Supp D$ is contained in the smooth part of $\mcC$, $D$ corresponds to a Cartier divisor on $\mcC$, and there is a section $s_D\in\Gamma(\mcC,\mcO_\mcC(D))$ defining $D$ and satisfying $s_D^2=F|_\mcC$. 
Put $\mcF:=\mcL|_\mcC\otimes\mcO_\mcC(-D)$ and let  $p_\mcF:\bm{L}_\mcF\to\mcC$ be the line bundle corresponding to $\mcF$. 
Note that $\mcF$ is of order $2$ and $\frac{t|_\mcC}{s_D}$ can be regarded as a section of $\Gamma(\bm{L}_\mcF,p_\mcF^\ast \mcF)$. 
We say that the unramified double cover $\varpi:\mcC_\mcF'\to\mcC$ given by
\[ \left( \frac{t|_\mcC}{s_D} \right)^2 - 1=0 \]
in $\bm{L}_\mcF$ and $\varpi:=p_\mcF|_{\mcC'}$ is induced by $\pi$. 
The morphism $\mcF\to\mcL|_\mcC$ given by multiplication of $s_D$ induces the morphism $\bm{L}_\mcF\to\bm{L}_\mcL|_\mcC$ over $\mcC$, which is given by multiplication of the value of $s_D$ to each fiber coordinate of $\bm{L}_\mcF$. 
This morphism induces $\mcC'_\mcF\to \pi^{-1}(\mcC)$, which is isomorphic over $\mcC\setminus\Supp D$. Hence we can deduce the structure of $\mcC'_\mcF$ by studying $\pi^{-1}(\mcC)$.
In the above cases,  explicit calculations of transition functions may be avoided when calculating the gluing data, which we see in the following example.

\begin{ex}\label{ex:2bitan}
Let $C_1, C_2\subset \PP^2$ be smooth conics intersecting transversely,  $T_1$, $T_2$, $T_3$, $T_4$ be the four bitangent lines to $C_1$, $C_2$ and $P_{ik}=T_i\cap C_k$ be the tangent points. Since $T_i+T_j$ has degree 2 and $(T_i+T_j)|_\mcC=2(P_{i1}+P_{i2}+P_{j1}+P_{j2})$, 
the ramified double cover $\pi_{ij}:S'_{ij}\rightarrow \PP^2$ of $\PP^2$ branched along $T_i+T_j$ induces an unramified double cover $\varpi_{ij}: \mcC'_{ij}\rightarrow \mcC$ of $\mcC=C_1+C_2$ as explained above by taking $\mcB=T_i+T_j$ and $D=P_{i1}+P_{i2}+P_{j1}+P_{j2}$. Note again that the  covers $\varpi_{ij}$ and $\pi_{ij}|_{\pi_{ij}^{-1}(\mcC)}$ of $\mcC$ are isomorphic outside the points $\{P_{i1}, P_{i2}, P_{j1}, P_{j2}\}$. 
Now let $S_{ij}\rightarrow S_{ij}^\prime$ be the canonical resolution of $\pi_{ij}$. Then $S_{ij}\cong \Sigma_2$, where $\Sigma_2$ is the Hirzebruch surface of degree 2, and we have the following diagram
\[
\begin{CD}
S'_{ij} @<<< S_{ij}\cong \Sigma_2 \\ 
@V{\pi_{ij}}VV   @VVV  \\ 
\PP^2@<<{\sigma_{ij}}< \widehat{\PP^2} \\ 
\end{CD}
\]
where $\sigma_{ij}$ is the blow up at the intersection point $T_i\cap T_j$. The pencil of lines through the intersection point induces the ruling of $\Sigma_2$. It can be readily checked that the preimages $C_1^\pm, C_2^\pm$ of $C_1, C_2$ in $S_{ij}$ are all linearly equivalent to $2F+\Delta_0$, since $C_i^++C_i^-\sim 4F+2\Delta_0$ and $C_i^{\pm}\cdot\Delta_0=0$,  where $F$ is the divisor class of fibers and $\Delta_0$ is the unique negative section with $\Delta_0^2=-2$. Hence $C_1^+\cdot C_2^+=C_1^+\cdot C_2^-=2$ and the splitting type of $(C_1, C_2; \mcB)$ is $(2, 2)$ since $C_i^{\pm}\cdot\Delta_0=0$. Also, since $\varpi_{ij}$ and $\pi_{ij}|_{\pi_{ij}^{-1}(\mcC)}$ are isomorphic outside $P_{ik}$, this implies that the gluing data of  $\varpi_{ij}:\mcC'_{ij}\rightarrow \mcC$ and  $\mcF_{ij}:=\mcO_{\PP^2}(1)|_\mcC\otimes \mcO_{\mcC}(-P_{i1}-P_{i2}-P_{j1}-P_{j2})$ is $(+,+,-,-)$ for a suitable choice of labels on the nodes $C_1\cap C_2$. Note again that changing all of the signs in  the gluing data gives isomorphic covers/sheaves, so the data  $(+, +, -, -)$ and $(-, -, +, +)$ give equivalent covers/sheaves.
Let $Q_1, Q_2, Q_3, Q_4$ be the nodes of $C_1+C_2$ and suppose that the gluing data of $\mcF_{12}$ is $(+, +, -, -)$ for $(Q_1, Q_2, Q_3, Q_4)$ in this order. Since $\mcF_{12}\otimes \mcF_{13}\cong \mcF_{23}$ and all of these sheaves  have two \lq$+$'s and two \lq$-$'s in the gluing data and are non-trivial, we can assume  that $\mcF_{12}, \mcF_{13}, \mcF_{23}$ are all distinct and the  gluing data of $\mcF_{13}$ is $(+, -, +, -)$ and the  gluing data of $\mcF_{23}$ is $(-, +, +, -)$, after changing the labels of $Q_3, Q_4$ if necessary. By the same argument, each triple $\mcF_{ij}, \mcF_{jk}, \mcF_{ik}$, $\{i, j, k\}\subset \{1, 2, 3, 4\}$ gives all three possible distinct invertible sheaves with two \lq$+$'s and two \lq$-$'s in the gluing data. Moreover $\mcF_{ij}\not\cong \mcF_{ik}$ if $j\not=k$.
Furthermore, this implies that $\mcF_{ij}\cong \mcF_{kl}$ for $\{i, j, k, l\}=\{1, 2, 3, 4\}$ as $\mcF_{kl}$ cannot be isomorphic to $\mcF_{ik}$ or $\mcF_{jk}$ and must be isomorphic to the remaining $\mcF_{ij}$ in the triple  $\mcF_{ij}, \mcF_{jk}, \mcF_{ik}$. 
\end{ex}

\begin{rem}
In Example~\ref{ex:2bitan},  $\mcO_{\mcC}(2P_{i1}+2P_{i2})\cong \mcL|_\mcC\cong \omega_\mcC$, where $\omega_\mcC$ is the dualizing sheaf of $\mcC$. Hence $\mcO_{\mcC}(P_{i1}+P_{i2})$ is a \emph{theta characteristics} of $\mcC$ (see \cite{harris82}). Now, $\mcF_{ij}=\mcL|_\mcC\otimes \mcO_{\mcC}(-P_{i1}-P_{i2}-P_{j1}-P_{j2})\cong \mcO_{\mcC}((P_{i1}+P_{i2})-(P_{j1}+P_{j2}))$ and is nothing but the difference between the odd theta characteristics $\mcO_{\mcC}(P_{i1}+P_{i2})$ and $\mcO_{\mcC}(P_{j1}+P_{j2})$.
\end{rem}

\begin{rem}
Since $\mcF_{ij}\cong \mcF_{kl}$, we have
$
\mcF_{ij}\otimes\mcF_{kl}\cong \mcO_\mcC,
$
which gives
\[
\mcO_\mcC\left(\sum_{i=1}^4(P_{i1}+P_{i2})\right)\cong \mcO_{\PP^2}(2)|_\mcC.\] This  implies Salmon's theorem. Namely the eight points of tangency $\{P_{11}, \ldots, P_{42}\}$ lie on a conic. See \cite[Theorem 3.3]{harris82} for details. Also another different proof can be found in \cite[Corollary~1.5]{masuya23}.
\end{rem}

\section{Unramified double covers of two conics induced by Poncelet transverses}\label{sec:poncelet-cover}

Let $C_1, C_2$ be smooth conics intersecting transversely as before. In this section, we consider the unramified double covers of $\mcC=C_1+C_2$ induced by $n$-sided Poncelet transverses. We note that it  is known that there exist $C_1, C_2$ intersecting transversely that admit an $n$-sided Poncelet transverse for any $n\geq 3$ (see \cite{barth-bauer}). 

\subsection{Degenerated Poncelet transverses}
\begin{figure}[h]
\begin{tabular}{cc}
         \begin{minipage}[h]{0.45\linewidth}
    \centering
    \includegraphics[width=5cm]{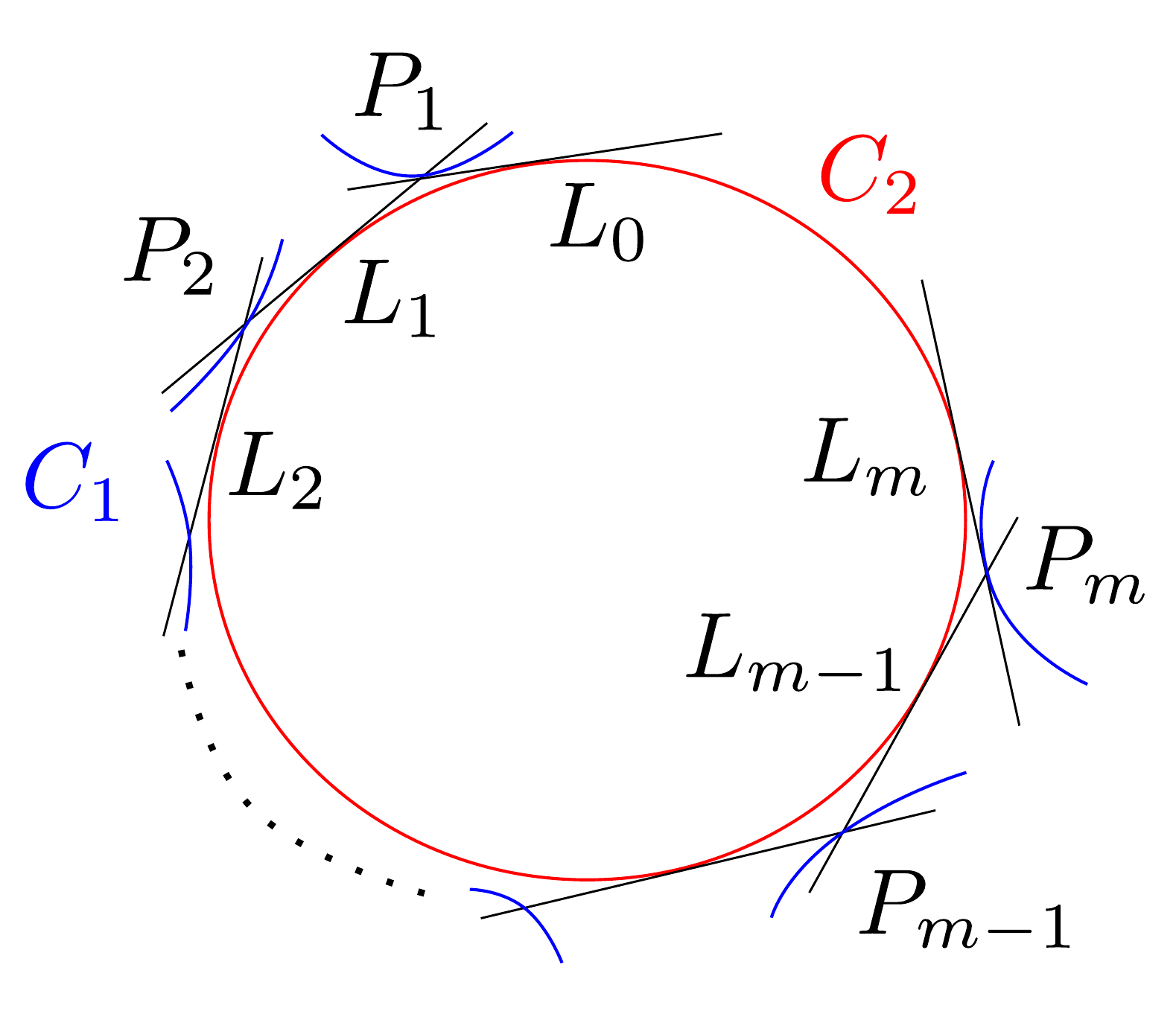}
    \subcaption{Case of $n=2m$. $L_0$ and $L_m$ are bitangent lines.}
  \end{minipage}&  
  \begin{minipage}[h]{0.45\linewidth}
    \centering
    \includegraphics[width=5.2cm]{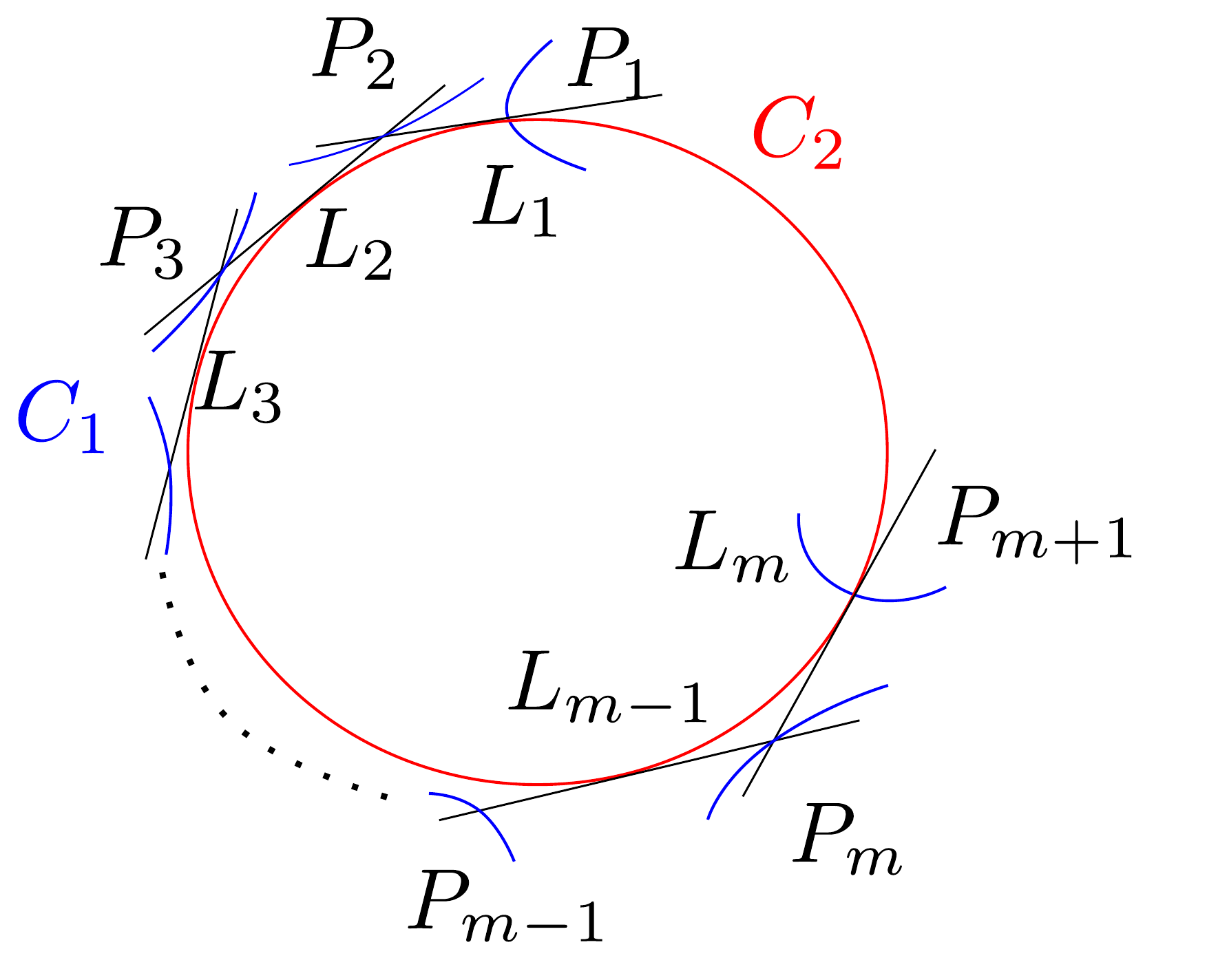}
    \subcaption{Case of $n=2m$. $L_1$ and $L_m$ are tangent to $C_2$ at points $P_1, P_{m+1}\in C_1 \cap C_2$.}
    \end{minipage} 
    \\
    \begin{minipage}[h]{0.45\linewidth}
    \centering
    \includegraphics[width=5cm]{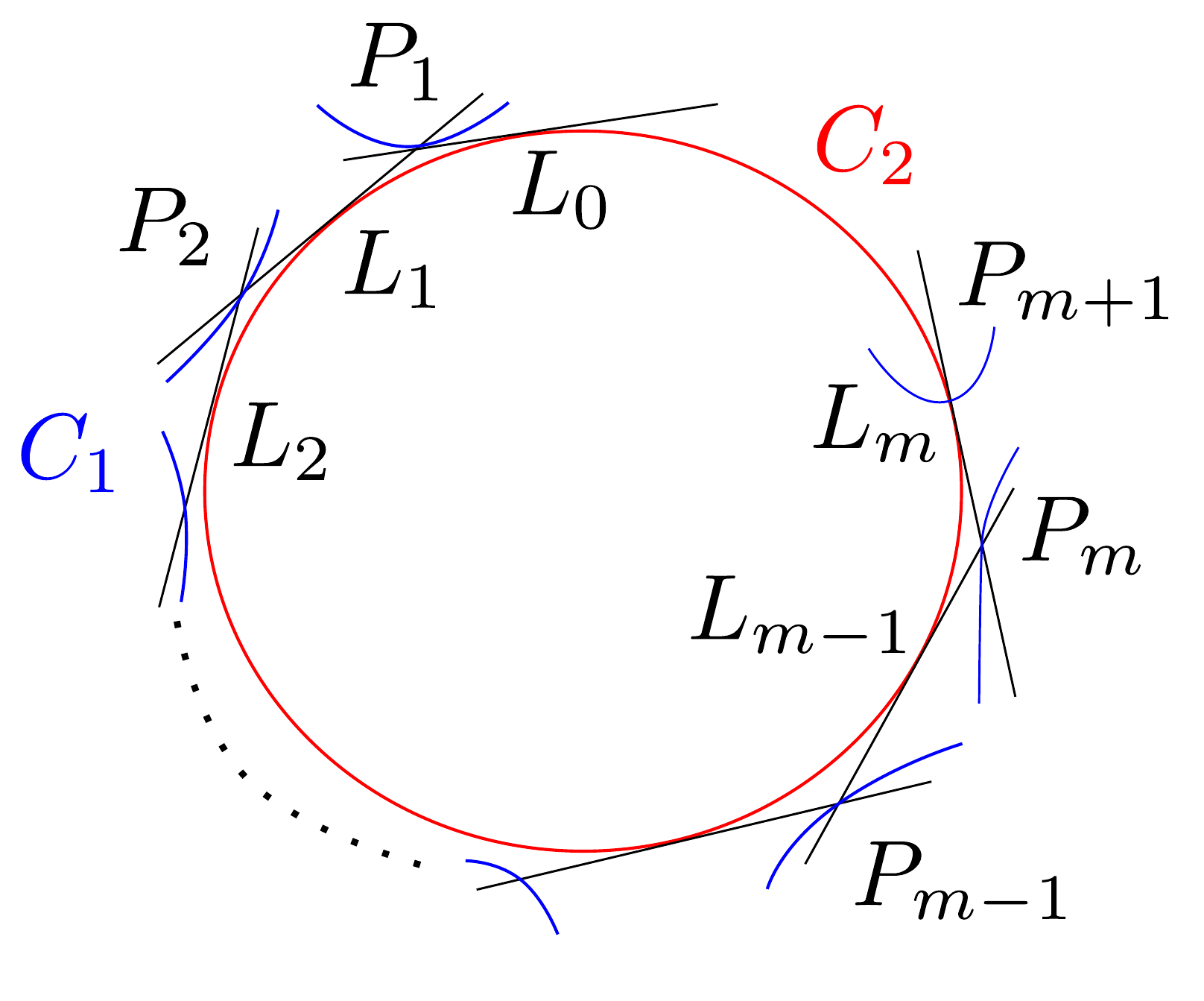}
    \subcaption{Case of $n=2m+1$. $L_0$ is a bitangent line and $L_m$ is a line tangent to $C_2$ at a point  $P_{m+1}\in C_1 \cap C_2$.}
\end{minipage}
\end{tabular}
\caption{Three types of degenerated $n$-sided Poncelet transverses.}
\label{fig:deg-trans}
\end{figure}

First, we study  degenerated $n$-sided Poncelet transverses. Let $C_1, C_2$ be smooth conics intersecting transversely with an $n$-sided Poncelet transverse $\sqcap_{n}: (P_1, L_1), \ldots, (P_n, L_n)$. If $\sqcap_n$ is degenerated,  there exists a pair $(P_i, L_i)$, $(P_j, L_j)$ ($i\not=j$) such that either $P_i=P_j$ or $L_i=L_j$. We can assume $i<j$ without loss of generality. Note that $(P_i, L_i)\not=(P_j, L_j)$ by the minimality of the period $n$. 
\begin{itemize}
    \item Suppose $P_i=P_j$. If  $P_i\in C_2$,  then there is only one unique line $L$ passing through $P_i$ and tangent to $C_2$. Then $L_i=L_j$ which contradicts the minimality of the period. Hence, we can assume $P_i\not\in C_2$ and that there exist two distinct tangent lines $L_i'$, $L_i''$, of $C_2$ passing through $P_i$. Since $L_i\not=L_j$ by the minimality of the period, we have  $\{L_i, L_j\}=\{L'_i, L_i''\}$. This implies that  $(P_i, L_i)$, $(P_j, L_j)$ are consecutive in the sequence and is of the form $(P_i, L_i)$, $(P_{i}, L_{i+1})$ and this can only occur if $L_i$ is a tangent line of $C_1$ and hence a bitangent line to $C_1+C_2$. 
    \item Suppose $L_i=L_j$. Then since $P_i$ and $P_j$ lie on the same line,  again we can assume that $(P_i, L_i)$, $(P_j, L_j)$ are consecutive in the sequence and is of the form $(P_i, L_i)$, $(P_{i+1}, L_i)$. If $P_{i+1} \not\in C_2$, then there exists two distinct lines through $P_{i+1}$ tangent to $C_2$ and $L_i\not=L_{i+1}$ which contradicts $L_i=L_j$. Hence $P_{i+1}\in C_2$ and $L_i=L_j$ is the unique tangent line of $C_2$ passing through $P_{i+1}$.
\end{itemize}
In both cases, the sequence is \lq reflected' at $(P_i, L_i)$, $(P_{i+1}, L_{i+1})$ and the points and lines leading up to this position appear in reverse order leading away. In order to be periodic, the sequence must be \lq reflected' once more to come back to $(P_i, L_i)$, $(P_{i+1}, L_{i+1})$. A \lq reflection' only occurs in the above two cases,  hence if $\sqcap_n$ is degenerated, then it must contain exactly two lines that are either a bitangent line or a  line tangent to $C_2$ at a point of $C_1\cap C_2$. A bitangent line will appear in the whole sequence only once and the other lines will appear exactly twice. Therefore we have the following:
\begin{itemize}
    \item If $n=2m$, two cases can occur. In the first case, $\sqcap_{2m}$ has two bitangent lines. The set of vertices consist of $m$ distinct points and the set of edges consist of two bitangents and $m-1$ general lines. In this case the transverse is of the form \[(P_1, L_1), \ldots, (P_m,L_m), (P_{m}, L_{m-1}),(P_{m-1}, L_{m-2}),\ldots, (P_{1}, L_{0})\]
    where $L_0$ and $L_m$ are the bitangent lines, under a suitable choice of labels (see Figure~\ref{fig:deg-trans} (a)). In the second case $\sqcap_{2m}$ has
    two lines each tangent to $C_2$ at a point of $C_1\cap C_2$. In this case the set of vertices consists of  $m+1$ distinct points and the set of edges consists of the two lines each tangent to $C_2$ at a point of $C_1\cap C_2$ and $m-2$ general lines. In this case the transverse is of the form 
    \[
    (P_1, L_1), (P_2, L_2),\ldots,  (P_m, L_m), (P_{m+1}, L_m), (P_m, L_{m-1}), \ldots, (P_3, L_2), (P_2, L_1) 
    \]
    where $P_1, P_{m+1}\in C_1\cap C_2$ and $L_1$ and $L_m$ are the lines tangent to $C_2$ at a point of $C_1\cap C_2$, under a suitable choice of labels (see Figure~\ref{fig:deg-trans} (b)).
    There exist two degenerated $2m$-sided transverses of each kind.
    \item If $n=2m+1$, $\sqcap_{2m+1}$ will have one bitangent line and one line tangent to $C_2$ at a point of $C_1\cap C_2$.
    The set of vertices consists of $m+1$ distinct points and the set of edges consist of the bitangent line, the line tangent to $C_2$ at a point in $C_1\cap C_2$ and $m-1$ general lines. The transverse is of the form 
    \[
    (P_1, L_1), \ldots,  (P_m, L_{m}), (P_{m+1}, L_{m}),(P_{m}, L_{m-1}) \ldots, (P_1, L_0)
    \]
    Where $L_0$ is the bitangent line, $P_{m+1}\in C_1\cap C_2$ and $L_{m}$ is the line tangent to $C_2$ at $P_{m+1}$, under a suitable choice of labels (see Figure~\ref{fig:deg-trans} (c)).
    There exist four degenerated $2m+1$-sided transverses of this kind.
\end{itemize}

%{\color{red} It would be nice to have some pictures here illustrating the degenerated Poncelet transverses.}

% \includegraphics[width=10cm]{n_sided_degenerated_3types.pdf}

% \begin{table}[h]
%     \centering
%     \begin{tabular}{|c|c|c|}
%         \hline
%         $n$ & \multicolumn{2}{|c|}{$C_1+C_2$ and a degenerated $n$-sided}   \\
%         \hline
%         \raisebox{6em}{$2m$} & 
%         $L_0$ and $L_m$ are bitangent lines. \includegraphics[width=5cm]{n_sided_degenerated_bitangents.pdf} &
%         \includegraphics[width=5cm]{n_sided_degenerated_intersections.pdf}\\
%         \hline
%         \raisebox{5em}{$2m+1$} & \multicolumn{2}{|c|}{\includegraphics[width=5cm]{n_sided_degenerated_odd_type.pdf}}\\
%         \hline
%     \end{tabular}
%     \caption{Three types of degenerated $n$-sideds}
%     \label{table:n_sided}
% \end{table}

\subsection{Deformation and degeneration of Poncelet transverses and line bundles of order two}\label{subsec:poncelet_cover}

Let $C_1, C_2$ be smooth conics intersecting transversely admitting a  $2m$-sided Poncelet transverse $\sqcap_{2m}$. In this subsection, we consider double covers of $\PP^2$ branched along  the lines of $\sqcap_{2m}$ and its relation with the induced unramified double covers of $\mcC=C_1+C_2$. We study the unramified double cover through a degeneration argument, where we deform  general $2m$-sided Poncelet transverses to a degenerated transverse with two bitangent lines. Note that we do not consider the other type of degeneration, as it will not induce an unramified double cover of $\mcC$.

Let $\sqcap_{2m}: (P_1, L_1), \ldots, (P_{2m}, L_{2m})$ be a general non-degenerated Poncelet transverse and let $\mcP_{2m}:=\sum_{i=1}^{2m}L_i$. Let $Q_i=C_2\cap L_i$ be the tangent points of $L_i$ and $C_2$ ($i=1, \ldots, 2m$). Let $\pi_{\mcP}: S'\rightarrow \PP^2$ be the double cover branched along $\mcP_{2m}$. Then, since $\mcP_{2m}$ has degree $2m$ and $\mcP_{2m}|_{\mcC}=2(\sum_{i=1}^{2m} P_i+\sum_{j=1}^{2m} Q_j)$,   $\pi_{\mcP}$ induces an unramified double cover  $\varpi_\mcP: \mcC_\mcP'\rightarrow \mcC$ as in Section~\ref{sec:linebdl2coverings}. The line bundle of order $2$ on $\mcC$ defining $\mcC'_\mcP$ is $\mcF_\mcP:=\mcO_{\PP^2}(m)|_\mcC\otimes \mcO_{\mcC}(-\sum_{i=1}^{2m} P_i-\sum_{j=1}^{2m} Q_j)$. We are interested in the structure of this double cover $\mcC'_\mcP$. By Theorem~\ref{thm:poncelet} (Poncelet's closure theorem), when we continuously move  $P_1$ on $C_1$ to a point $P'$ that is a tangent point of a bitangent line of $\mcC=C_1+C_2$, the  $2m$-sided Poncelet transverse originating at $P_1$ continuously deforms along with $P_1$ to the degenerated Poncelet transverse originating at $P'$. Let $P'_1=P'$ and let 
\[(P'_1, L'_1), \ldots, (P'_m,L'_m), (P'_{m}, L'_{m-1}),(P'_{m-1}, L'_{m-2}),\ldots, (P'_{1}, L'_{0})\]
be the degenerated transverse originating at $P_1'$ where $L_0'$ and $L_{m}'$ are bitangent lines of $\mcC=C_1+C_2$. Let $Q_i'=L'_i\cap C_2$ ($i=0, \ldots, m$) be the tangent points of $L'_i$ and $C_2$.
The correspondence between the lines, vertices and tangent points under the degeneration are as follows:
\begin{align*}
 P_i, P_{2m+1-i} &\rightarrow P'_i \quad (i=1,\ldots, m),\\
L_i, L_{2m-i} &\rightarrow L'_i \quad  (i=1, \ldots, m-1), \quad L_m \rightarrow L'_m, \quad L_{2m} \rightarrow L'_0\\
Q_i, Q_{2m-i} &\rightarrow Q'_i \quad  (i=1, \ldots, m-1), \quad Q_m \rightarrow Q'_m, \quad Q_{2m} \rightarrow Q'_0
\end{align*}

%{\color{red} Make a nice table or figure or something describing the correspondence when degenerating}

%{\color{red} Do we need any remark saying that we only need to consider
%{{\it bitangent} degenerations?}
\begin{table}[h]
    \centering
    \begin{tabular}{|c|c|c|}
       % \hline
        % & \multicolumn{2}{|c|}{$C_1+C_2$ and a $6$-sided Poncelet transverse}  \\
         \hline
         \raisebox{3.5em}{Non-degenerated}& 
        \includegraphics[width=4cm]{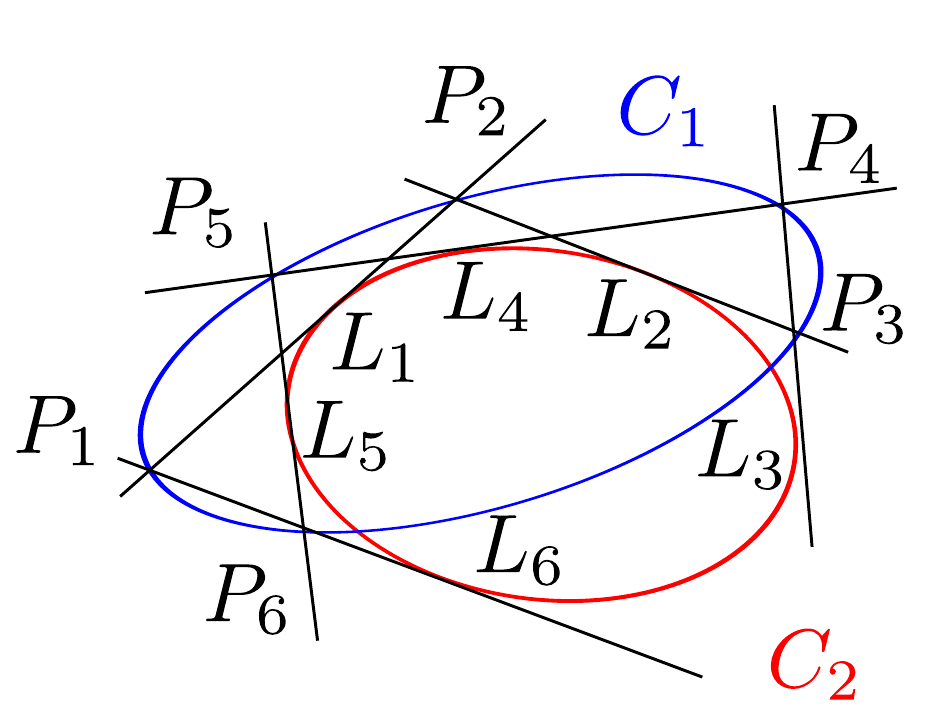} & 
         \raisebox{3.5em}{$\begin{array}{c}
            (P_1,L_1) \to (P_2,L_2)\to (P_3,L_3)\\
            \to (P_4,L_4) \to (P_5,L_5) \to (P_6,L_6) 
        \end{array}$} \\
        \hline
        \raisebox{3.5em}{Degenerated}& 
        \includegraphics[width=4.5cm]{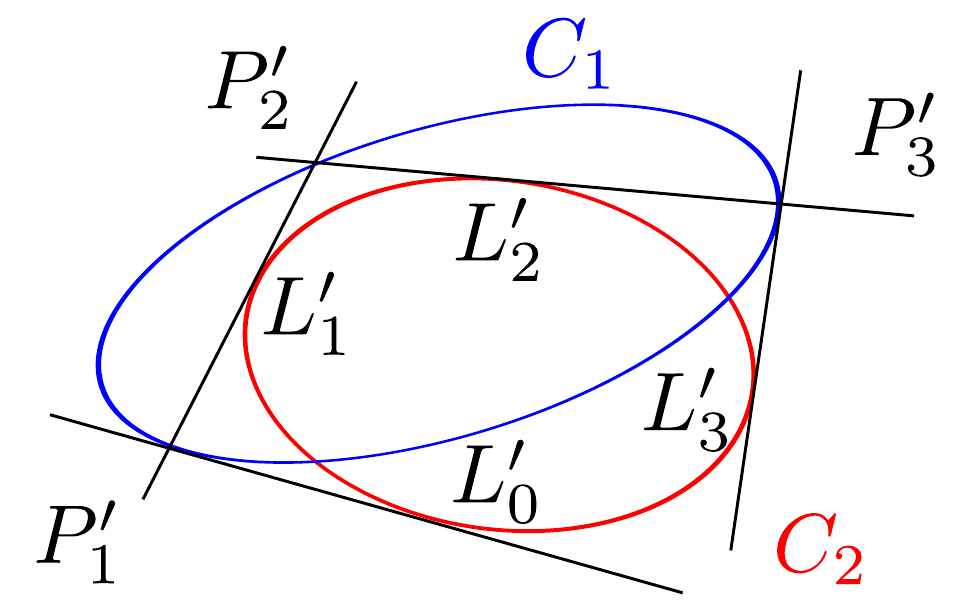} & \raisebox{3.5em}{$\begin{array}{c}
            (P_1',L_1') \to (P_2',L_2') \to(P_3',L_3')\\
            \to (P_3',L_2') \to (P_2',L_1') \to (P_1',L_0') 
        \end{array}$} \\
        \hline
    \end{tabular}
    % \caption{Case of $m=3$. A non-degenerated $6$-sided and a degenerated $6$-sided with two bitangent lines.}
    % \caption{A non-degenerated $6$-sided Poncelet traverse and a degenerated $6$-sided Poncelet traverse with two bitangent lines.}
     \caption{$C_1+C_2$ and the degeneration of a $6$-sided Poncelet transverse.}
    \label{table:6_sided}
\end{table}

%$1\leq i\leq m+1$,  the lines $L_i$ deform to  $L'_i$ and for $m+2\leq i \leq 2n$, the lines $L_i$ deform to $L'_{2m+2-i}$. 

Here the points  $P_1, \ldots, P_{2m}$ and $Q_1, \ldots, Q_{2m}$ are continuously deformed on the smooth part of $\mcC$, while preserving the condition that $\mcO_{\PP^2}(m)|_\mcC\otimes \mcO_{\mcC}(-\sum_{i=1}^{2m} P_i-\sum_{j=1}^{2m} Q_j)$ is of order $2$. On the other hand, the set of invertible sheaves of order $2$ of $\mcC$ is  isomorphic to $(\ZZ/2\ZZ)^{\oplus 3}$ (see the discussions in Section~\ref{sec:linebdl2coverings} or \cite[Section~3, 3a]{harris82}) and is finite and discrete. By Lemma~\ref{lem:rigidity-of-top-2cover} the gluing data  of the induced unramified double cover is constant under the deformation, hence by Lemma~\ref{lem:rigidity_of_sheaves} the gluing data of the invertible sheaves of   order 2 must also be constant under the deformation and  the isomorphism classes of the sheaves $\mcF_{\mcP}$  must be constant. Hence, we have 
\[ \mcO_{\PP^2}(m)|_\mcC\otimes \mcO_{\mcC}\left(-\sum_{i=1}^{2m} P_i-\sum_{j=1}^{2m} Q_j\right)\cong \mcO_{\PP^2}(m)|_\mcC\otimes \mcO_{\mcC}\left(-\sum _{i=1}^m 2P'_i-\sum_{j=1}^{m-1} 2Q'_j-Q'_0-Q'_{m}\right). \] 
Furthermore, since $(L_1'+\cdots+L_{m-1}')|_\mcC=P_1'+\sum_{i=2}^{m-1}2P_i'+P_m'+\sum_{j=1}^{m-1} 2Q_j'$ and $\mcO_{\PP_2}(m-1)|_{\mcC}\cong \mcO_{\mcC}(P_1'+\sum_{i=2}^{m-1}2P_i'+P_m'+\sum_{j=1}^{m-1} 2Q_j')$, we have 
\[ \mcO_{\PP^2}(m)|_\mcC\otimes \mcO_{\mcC}\left(-2\sum _{i=1}^m P'_i-\sum_{j=1}^{m-1} 2Q'_j-Q'_0-Q'_{m}\right) \cong \mcO_{\PP^2}(1)|_\mcC\otimes \mcO_{\mcC}(- P'_1-P'_m-Q'_0-Q'_{m}). \]
The points $ P'_1, Q'_0$ are the tangent points of the bitangent line $L_0'$ and  the points $P'_m,Q'_{m}$ are the tangent points of the bitangent line $L_{m}'$. Hence we see that the structure of 
$\mcF_\mcP$ and the associated unramified double cover $\varpi_\mcP:\mcC'_\mcP\rightarrow\mcC$ induced by $\pi_\mcP$ is identical to that of the double cover in Example~\ref{ex:2bitan} associated to the bitangent lines $L_0'+L_{m}'$. Summing up these arguments, we have the following Lemma, where $\varpi_{ij}$ is the unramified  double cover of $\mcC$ induced by  the double cover of $\PP^2$ branched along $T_i+T_j$ as  defined in Example~\ref{ex:2bitan}.
\begin{lem}\label{lem:trans_bitan}
Under the above settings and notation, let $T_1, T_2, T_3, T_4$ be the bitangent lines to $C_1, C_2$ labeled so that the pairs $T_1, T_2$ and $T_3, T_4$ each lie in the same degenerated $2m$-sided transverse. Then the  unramified double covers $\varpi_{\mcP}$, $\varpi_{12}$, $\varpi_{34}$ of $\mcC=C_1+C_2$ are all isomorphic.
\end{lem}

\begin{rem}
The isomorphism between $\varpi_{12}$ and $\varpi_{34}$ has already been observed in Example~\ref{ex:2bitan}, regardless of the existence of a $2m$-sided Poncelet transverse.
\end{rem}

\section{Proof of Main Theorem}\label{sec:proof-main}
In this section, we  prove Theorem~\ref{thm:main}. Let  $C_1, C_2$ be smooth conics intersecting transversely that admits a $2m$-sided Poncelet transverse. Let $T_1, T_2, T_3, T_4$ be bitangent lines of $\mcC=C_1+C_2$ labeled so that the pairs $T_1$, $T_2$ and $T_3$, $T_4$ each lie in the same degenerated $2m$-sided transverse. Let, $\sqcap_{2m}: (P_1, L_1), \ldots, (P_{2m}, L_{2m})$ be a non-degenerated transverse, $\mcP_{2m}:=\sum_{i=1}^{2m}L_i$ and let
\[
\mcC_{ij}:=C_1+C_2+\mcP_{2m}+T_i+T_j \quad (\{i, j\}\subset \{1, 2, 3, 4\})
\]
as in the Introduction.

\begin{lem}\label{lem:comb}
The combinatorial types $\comb(\mcC_{ij})$ of $\mcC_{ij}$ are the same for all $\{i, j\}\subset \{1, 2, 3, 4\}$ and any choice of non-degenerated transverse $\sqcap_{2m}$.
\end{lem}

\begin{proof} Let $\sqcap_{2m}: (P_1, L_1), \ldots, (P_{2m}, L_{2m})$ be a non-degenerated transverse. Since all of the lines $L_1, \ldots, L_{2m}$ and $T_1, T_2, T_3, T_4$ are tangent lines of $\mcC_2$, no three are concurrent.  A line $L_i$ and a bitangent $T_j$ cannot intersect on $C_1$ as we have assumed that $L_i$ lies in a non-degenerated transverse and $T_j$ lies in a degenerated transverse. They cannot intersect on $C_2$ as well because they are distinct tangent lines of $C_2$. Hence the combinatoral types are the same. \end{proof}

Let $\mcB_{ij}:= \mcP_{2m}+T_i+T_j \quad (\{i, j\}\subset\{1, 2, 3, 4\})$
and let
$\pi_{\mcB_{ij}}:S'\rightarrow \PP^2$ be the double cover of $\PP^2$ branched along $\mcB_{ij}$. 
\begin{lem}\label{lem:splitting_type}
Under the labeling above, 
the splitting type of $(C_1, C_2; \mcB_{ij})$ is $(0,4)$ if $\{i, j\}=\{1, 2\}$ or $\{3, 4\}$ and is $(2, 2)$ otherwise.  
\end{lem}
\begin{proof} Since $\mcB_{ij}$ can be viewed as a sum of $\mcP_{2m}$ and $(T_i+T_j)$, by the discussions in Section~\ref{sec:linebdl2coverings}, the cover $\pi_{\mcB_{ij}}$ induces an unramified double cover of $\mcC=C_1+C_2$, whose structure is given by the product of the covers $\varpi_{\mcP}$ of Section~\ref{subsec:poncelet_cover} and $\varpi_{ij}$ of Example~\ref{ex:2bitan}. Then since  $\varpi_{\mcP}$, $\varpi_{12}$, $\varpi_{34}$ are isomorphic by Lemma~\ref{lem:trans_bitan}, $\mcB_{ij}$ induces the trivial unramified double cover if $\{i, j\}=\{1, 2\}$ or $\{3, 4\}$, and otherwise induces a non-trivial unramified double cover with gluing data $(+, +, -, -)$ for a suitable choice of labels of the nodes. 
\end{proof}
Now,  Lemma~\ref{lem:comb} and \ref{lem:splitting_type} together with Proposition~\ref{prop:splitting-type} give the proof of Theorem~\ref{thm:main}.

\bibliographystyle{plain}
\bibliography{biblio.bib}

\noindent Shinzo BANNAI \\
Department of Applied Mathematics, Faculty of Science, \\
Okayama University of Science, 1-1 Ridai-cho, Kita-ku,
Okayama 700-0005 JAPAN \\
{\tt bannai@ous.ac.jp}\\

\noindent  Ryosuke MASUYA, Hiro-o TOKUNAGA and Emiko YORISAKI\\
Department of Mathematical  Sciences, Graduate School of Science, \\
Tokyo Metropolitan University, 1-1 Minami-Ohsawa, Hachiohji 192-0397 JAPAN \\
{\tt tokunaga@tmu.ac.jp}\\

\noindent Taketo SHIRANE\\
Department of Mathematical Science, Faculty of Science and Technology, \\
Tokushima University, 2-1 Minamijyousanjima-cho, Tokushima 770-8506, JAPAN\\
{\tt shirane@tokushima-u.ac.jp}

\end{document}